\documentclass[a4paper,11pt]{article}
\usepackage{textcomp}
\usepackage{amsfonts}
\usepackage[dvips]{graphicx}
\usepackage{amsmath}
\usepackage{amssymb}
\usepackage{mathrsfs}

\begin{document}

\title{\textbf{Automorphic Forms and Reeb-Like Foliations on Three-Manifolds}}

\author{
Yi SONG, Xu XU and Stephen P. BANKS\\
Dept. of Automatic Control \& Systems Engineering, \\
University of Sheffield,\\
Mappin Street, Sheffield, \\
S1 3JD.\\
email: s.banks@sheffield.ac.uk \\
}

\maketitle

\newtheorem{theorem}{Theorem}[section]
\newtheorem{lemma}{Lemma}[section]
\newtheorem{definition}{Definition}[section]
\newtheorem{example}{Example}[section]

\thispagestyle{empty}

\begin{abstract}
In this paper, we consider different ways of generating dynamical
systems on 3-manifolds. We first derive explicit differential
equations for dynamical systems defined on generic hyperbolic
3-manifolds by using automorphic function theory to uniformize the
upper half-space model. It is achieved via the modification of the
standard \textit{Poincar\'{e}} theta series to generate systems
invariant within each individual fundamental region such that the
solution trajectories match up on the appropriate sides after the
identifications which generate a hyperbolic 3-manifold. Then we
consider the gluing pattern in the conformal ball model. At the end
we shall study the construction of dynamical systems by using
the Reeb foliation.
\end{abstract}

\section{Introduction}

Nonlinear dynamical systems are defined globally on manifolds.
Consider a system

\begin{equation}\label{general system}
\dot{x}=Ax
\end{equation}

\noindent where $x\in \mathbb{R}^n$. The phase-space portrait of
this system is defined as all the solution curves in $\mathbb{R}^n$
(see, e.g., [Perko, 1991].) Geometrically, these curves determine
the motion of all the points in the space under this specific
dynamical system.

\begin{figure}[!h]
\centering
\includegraphics[width=2in]{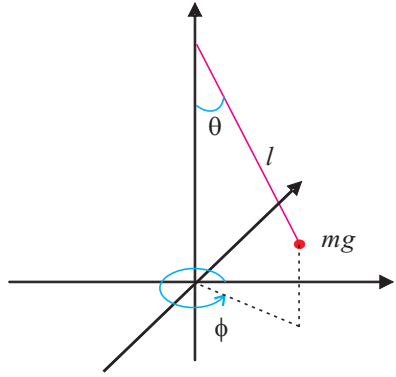}
\caption{A spherical pendulum}\label{spherical pendulum}
\end{figure}

Moreover, the global manifold on which a system is stratified can be
reached by studying the phase-space portrait. For example, given a
spherical pendulum (see fig.~\ref{spherical pendulum} for
illustration), it has two degrees of freedom, and the
\textit{Lagrangian} for this system is

\begin{equation}
L=\frac{1}{2}m\big(l^2{\dot{\theta}}^2+l^2(\sin{\theta})^2{\dot{\phi}}^2
\big)+mgl\cos{\theta}.
\end{equation}

\noindent The \textit{Euler-Lagrange} equations give

\begin{equation}
\left\{
\begin{array}{l}
\frac{d}{dt}(ml^2\dot{\theta})-ml^2\sin{\theta}\cos{\theta}{\dot{\phi}}^2+mgl\sin{\theta}=0 \\
\frac{d}{dt}\big(ml^2({\sin{\theta}})^2\dot{\phi}\big)=0
\end{array}\right.
\end{equation}

\noindent so the system is given by two equations of motion, i.e.,

\begin{equation}
\left\{
\begin{array}{l}
\ddot{\theta}={\dot{\phi}}^2\sin{\theta}\cos{\theta}-g/l\cdot
\sin{\theta} \\
\ddot{\phi}=-({2\cos{\theta} \cdot \dot{\theta} \cdot
\dot{\phi}})/{\sin{\theta}}
\end{array}\right.
\end{equation}

In phase-space coordinates, if set $x_1=\theta$,
$x_2=\dot{\theta}=\omega_{\theta}$, $x_3=\phi$,
$x_4=\dot{\phi}=\omega_{\phi}$, we then have

\begin{equation}
\left\{
\begin{array}{l}
\dot{x}_1=x_2 \\
\dot{x}_2={x_4}^2 \sin{x_1} \cos{x_1} - g/l \cdot \sin{x_1} \\
\dot{x}_3=x_4 \\
\dot{x}_4=-({2 \cos{x_1} \cdot x_2 \cdot x_4})/{\sin{x_1}}
\end{array}\right.
\end{equation}

\begin{figure}[!h]
\centering
\includegraphics[width=2in]{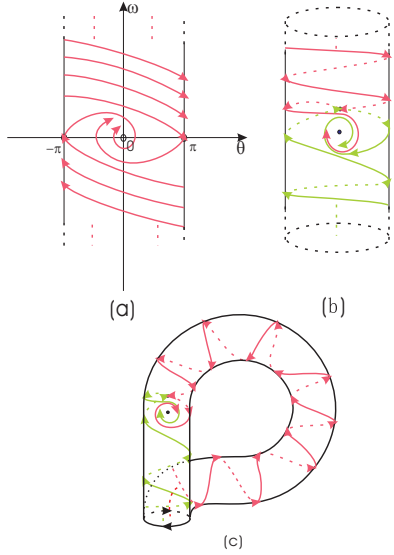}\\
\caption{A single pendulum is defined on a \textit{Klein
bottle}}\label{simple pendulum}
\end{figure}

\noindent This is a 4-dimensional system. Now assume $x_3={\phi}=k$,
where \textit{k} is a constant, consequently $\dot{x}_3=0$ and
${\dot{x}}_4={\dot{\omega}}_{\phi}=0$, the system will then become

\begin{equation}
\left\{
\begin{array}{l}
\dot{x}_1=x_2 \\
\dot{x}_2=- g/l \cdot \sin{x_1} \\
\dot{x}_3=0 \\
\dot{x}_4=0
\end{array}\right.
\end{equation}

\noindent which stands for a single pendulum that sits on a
2-dimensional plane. It is known that this system is defined on a
\textit{Klein bottle}, (see [Banks \& Song, 2006] and
fig.~\ref{simple pendulum} for an illustration.)

\begin{figure}[!h]
\centering
\includegraphics[width=3in]{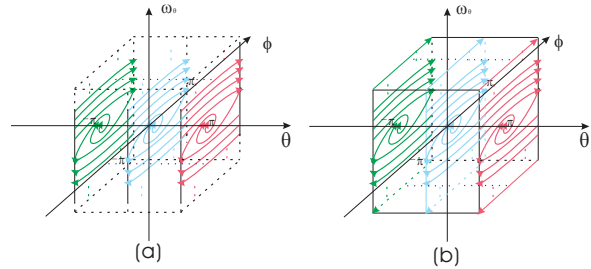}\\
\caption{Phase-plane portrait of the spherical pendulum when
$\omega_{\phi}=k$}\label{special spherical pendulum}
\end{figure}

If set ${\dot{x}}_3=\omega_{\phi}=k$, i.e., the system has a fixed
angular velocity in $\phi$, the system then becomes

\begin{equation}\label{order reduced}
\left\{
\begin{array}{l}
\dot{x}_1=x_2 \\
\dot{x}_2={x_4}^2 \sin{x_1} \cos{x_1} - g/l \cdot \sin{x_1} \\
\dot{x}_3=k \\
\dot{x}_4=0
\end{array}\right.
\end{equation}

\noindent which is essentially a 3-dimensional hyperplane given by
$x_4=k$ within the 4-dimensional space. Moreover, since the vector
field is periodic in both $x_1$ and $x_3$ with period $2\pi$, it is
naturally defined within the cube
\begin{eqnarray*}
C:\{(x_1,x_2,x_3): -\pi \le x_1 \le \pi, -\infty < x_2 < \infty,\\
-\pi \le x_3 \le \pi\}
\end{eqnarray*}

\noindent as shown in fig.~\ref{special spherical pendulum}(a).

Note that the phase-space portrait is a 2-dim single pendulum that
sits on different slices defined by $\phi=k$, and because  we know
that $\theta, \phi=\pi$ and $\theta, \phi=-\pi$ are physically the
same respectively, we can identify them by pairing the opposite
sides via translation.

\begin{figure}[!h]
\centering
\includegraphics[width=3in]{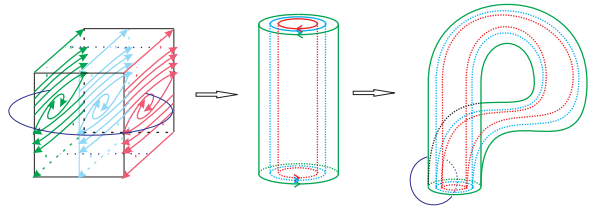}\\
\caption{Construction of a 3-dimensional solid \textit{Klein
bottle}}\label{solid klein bottle}
\end{figure}

In order to define the system on a compact manifold, we compress the
infinite cube to a finite one, as shown in fig.~\ref{special
spherical pendulum}(b), and since the dynamics at the two ends are
pointing the opposite directions, the identification will result in
a self-intersection in the 3-dimensional Euclidean space,
fig.~\ref{solid klein bottle} shows an embedding in $\mathbb{R}^3$.

Thus, we obtain the 3-manifold on which this special spherical
pendulum is defined. We call it 3-dimensional \textit{solid Klein
bottle}.

In [Banks \& Song, 2006], we showed that a dynamical system on a
two-dimensional surface is given by a generalized automorphic
function \textit{F}. In this paper, we will extend the previous
result and propose to show how to generalise explicit differential
equations that naturally have global behaviour on 3-manifolds. Again
we will use the theory of automorphic functions to achieve it.

\section{Geometric 3-Manifolds}

We shall now give a brief resum\'e of 3-manifolds which will be
needed in the following sections. Note that all the results are
well-known, for example in [Ratcliffe, 1994].

\begin{definition}
A $3$-manifold \textit{M} without boundary is a $3$-dimensional
Hausdorff space that is locally homeomorphic to $\mathbb{E}^3$,i.e.,
for every point $x$ $(x\in M)$ there exists a homeomorphism that
maps a neighbourhood \textit{A} of \textit{x} onto the
$3$-dimensional Euclidean space; while if \textit{M} has a boundary,
then the homeomorphism maps \textit{A} onto the upper-half
$3$-dimensional Euclidean space ~$\overline{\mathbb{U}^3} = \{x\in
\mathbb{E}^3: x_3\ge 0 \}$.
\end{definition}

Equivalently, a 3-manifold \textit{M} is called a geometric 3-space.
Assume that $\Gamma$ is a group which acts on a 3-dimensional
geometric space \textit{X}, then

\begin{definition}
The orbit space of the action $\Gamma$ on $X$ is the set of
$\Gamma$--orbits,
\[
X/\Gamma=\{\Gamma_X:x\in \Gamma\},
\]
with the metric topology being the quotient topology, and the
quotient map given by

\[
\pi : X\to X/\Gamma.
\]
\end{definition}

Moreover, if $\Gamma$ is a discrete group of isometries of
\textit{X}, then $\Gamma$ is discontinuous and called a
3-dimensional \textit{Fuchsian} group. In fact, it defines a
fundamental region \textit{F} of \textit{X} which, together with its
congruent counterparts, generates a tessellation of \textit{X}.

\begin{definition}
For a discrete group $\Gamma$ of isometries of a geometric space
\textit{X}, a subset \textit{F} of \textit{X} is a fundamental
region if and only if
\begin{enumerate}
\item the set \textit{F} is open in \textit{X};
\item the members of $\{ gF:g\in \Gamma\}$ are mutually disjoint;
\item $X=\cup\{g\overline{R}: g \in \Gamma \}$
\end{enumerate}
\end{definition}

\begin{figure}[!h]
\centering
\includegraphics[width=2in]{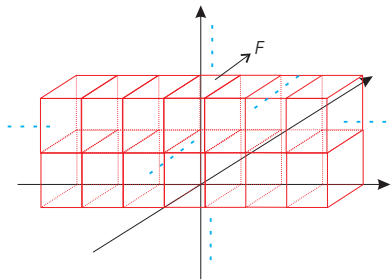}
\caption{Tessellation of $\mathbb{E}^3$ by $\Gamma_a$ generated via
translation $e_i$}\label{cube}
\end{figure}

For example, let $\tau_i$ be the translation of $\mathbb{E}^3$ by
$e_i$ for $i=1,2,3$, then $\{\tau_1,\tau_2,\tau_3\}$ defines a
discrete subgroup $\Gamma_a$ of $I(\mathbb{E}^3)$. A fundamental
region for $\Gamma_a$ will be the open unit cube in $\mathbb{E}^3$,
as shown in fig.~\ref{cube}, in fact, $\Gamma_a$ generates a
tessellation of $\mathbb{E}^3$.

If $\Gamma$ acts \textit{freely} on \textit{X}, the orbit space
$X/\Gamma$ is then a 3-manifold which can also be called an
\textit{X-space-form}. Also, by assuming \textit{G} is a group of
similarities of a 3-dimensional geometric space \textit{X} and
\textit{M} is a 3-manifold, we have

\begin{definition}
An $(X,G)$--atlas for \textit{M} is a group of maps
\[
\Phi=\{\phi_i:U_i\to X \}_{i\in \mathcal{I}}
\]

\noindent such that:

\begin{enumerate}
\item The set $U_i$ is an open connected subset of \textit{M} for
each \textit{i}.
\item $\phi_i$ maps $U_i$ homeomorphically onto an open subset of
\textit{X} for each \textit{i}.
\item $\bigcup_{i\in \mathcal{I}}U_i=M$
\item If $U_i$ and $U_j$ overlap, then the map
\[
\phi_j \phi_i^{-1}: \phi(U_i\cap U_j) \rightarrow \phi_j (U_i \cap
U_j),
\]
\noindent agrees in a neighbourhood of each point of its domain with
an element in \textit{G}.
\end{enumerate}
\end{definition}

Note that $\Phi$ consists of the \textit{charts} of the
\textit{(X,G)--atlas}. An \textit{(X,G)--structure} is then defined
as the maximal \textit{(X,G)--atlas} for \textit{M}. Hence a
3-manifold \textit{M} with an \textit{(X,G)--structure} is called an
\textit{(X,G)--manifold}. It is well-known (e.g., [Ratcliffe, 1994])
that the orbit space $X/\Gamma$, together with the induced
\textit{$(X,\Gamma)$--atlas}, is an \textit{$(X,\Gamma)$-manifold}.
Furthermore, we can obtain this 3-manifold by gluing one fundamental
region \textit{F} along the corresponding sides.

Let $\mathcal{F}$ be a family of fundamental regions in a geometric
space \textit{X} and $\Gamma$ be a group of isometries of
\textit{X}. We can then construct the $(X,\Gamma)$--manifold by
applying the \textit{$\Gamma$-side-pairing}.

\begin{definition}
A \textit{$\Gamma$-side-pairing} for $\mathcal{F}$ is a subset of~
$\Gamma$,

\[
\Gamma=\{\tau_s: \textrm{S is a side of one fundamental region in}~
\mathcal{F}\},
\]

\noindent such that for each side \textit{S} in $\mathcal{F}$,

\begin{enumerate}
\item there exists a side $S'$ in $\mathcal{F}$ that satisfies
$\tau_s(S')=S$,
\item $\tau_{s'}=\tau_s^{-1},$
\item if \textit{S} is a side of \textit{F} in $\mathcal{F}$ and ${S'}$ is a side of
${F'}$, then
\[
F\cap g_s(F')=S.
\]
\end{enumerate}
\end{definition}

\begin{figure}[!h]
\centering
\includegraphics[width=2in]{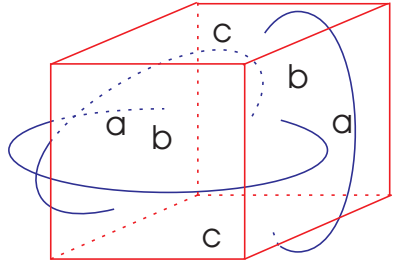}
\caption{Construction of the cubical Euclidean
3-torus}\label{3-torus}
\end{figure}

The elements of $\Gamma$ are called the side-pairing transformations
of $\mathcal{F}$, which generates an equivalence relation on the set
$\prod = \bigcup_{F\in \mathcal{F}}F$, i.e., the cycles of $\Gamma$.
Moreover, $S'$ is uniquely determined by $S$. So if the
$\Gamma$-side-pairing is proper, i.e., each cycle of $\Gamma$ is
finite and has solid angle sum $4\pi$, then by choosing two
fundamental regions in $\mathcal{F}$, say $F$ and $F'$, the elements
in $\Gamma$ will associate each side in $F$ with a unique one in
$F'$, identifying the corresponding sides together will eventually
generate a 3-manifold with an \textit{$(X,\Gamma)$-structure}
attached. For instance, as in the previous example, after pairing
the opposite sides of the unit cube by translations $\Gamma$, we
effectively end up with a 3-manifold \textit{M} which is known as
the \textit{cubical Euclidean 3-torus} (see fig.~\ref{3-torus} for
illustration).

\section{Automorphic Functions and Systems on Hyperbolic 3-Manifolds}

In this section we shall first give a brief resum\'e of automorphic
functions. More details can be found in, for example, [Ford, 1929;
Ratcliffe, 1994].

To denote the points in $\mathbb{R}^3$, we use the following
coordinates:

\begin{eqnarray*}
\mathbb{R}^3 &=& \mathbb{C}\times(-\infty,\infty) \\
 &=& \{(z,r)\mid z\in \mathbb{C}, r\in \mathbb{R}\}\\
 &=& \{(x,y,r)\mid x,y,r\in \mathbb{R}\}
\end{eqnarray*}

Also, we can think $\mathbb{R}^3$ as a subset of Hamilton's
quaternions $\mathcal{H}$, so a point \textit{p} $(p\in
\mathbb{R}^3)$ can be expressed as a quaternion whose fourth term
equals to zero, i.e.,

\[
p=(z,r)=(x,y,r)=z+rj,
\]

\noindent where $z=x+yi$ and $j=(0,0,1)$, then

\begin{definition}
A M\"obius transformation of ~${\widehat{\mathbb{R}}}^3$ is a finite
composition of reflections of ~${\widehat{\mathbb{R}}}^3$ in
spheres, where ~${\widehat{\mathbb{R}}}^3$ is the one-point
compactification of ~$\mathbb{R}^3$, i.e.,
\[
{\widehat{\mathbb{R}}}^3=\mathbb{R}^3\cup \{\infty\}.
\]
\end{definition}

It is exactly the linear fractional transformations of the form
\begin{equation}\label{mobius transformation}
T=\frac{ap+b}{cp+d}.
\end{equation}

\noindent where $a,b,c,d\in \mathbb{R}^3$ and $ad-bc\ne 0$.

A \textit{M\"obius} transformation is a conformal map of the
extended 3-space, (i.e., Riemann 3-manifold), denoted by
$\textbf{Aut}(\widehat{\mathbb{R}}^3)$. Moreover, (\ref{mobius
transformation}) can be represented in terms of a matrix

\begin{equation}
G=\left(
     \begin{array}{cc}
       a & b \\
       c & d \\
     \end{array}
   \right).
\end{equation}

\noindent In fact there exists a group homeomorphism:
$\textrm{GL}(2,\mathbb{R}^3) \to
\textbf{Aut}({\widehat{\mathbb{R}}}^3)$ given by

\[
\left(
     \begin{array}{cc}
       a & b \\
       c & d \\
     \end{array}
   \right) \longrightarrow T,
\]

\noindent which becomes an isomorphism on the projective special
linear group $\textrm{PSL}(2,\mathbb{R}^3)$ \big(i.e., those
elements of $\textrm{GL}(2,\mathbb{R}^3)$ of positive determinant
modulo the scalar matrices\big).

It is known that 3-dimensional hyperbolic space (or 3-dimensional
hyperbolic manifold) is the unique 3-dimensional simply connected
Riemannian manifold with constant sectional curvature $-1$ (see,
e.g., [Elstrodt, 1998]). Also, since \textit{M\"obius}
transformations are defined on Riemannian manifold, they can be used
to generate a discrete group of discontinuous isometries, $\Gamma$,
of the upper half-space $\mathbb{U}^3$, where

\begin{eqnarray*}
\mathbb{U}^3: &=& \mathbb{C}\times (0,\infty) \\
 &=& \{(z,r)\mid z\in \mathbb{C}, r>0\} \\
 &=& \{(x,y,r)\mid x,y,r\in \mathbb{R},r>0\}.
\end{eqnarray*}

Note that $\mathbb{U}^3$ is a model for hyperbolic space, so we can
use $\Gamma$ to tessellate $\mathbb{U}^3$ and obtain a 3-manifold,
\textit{M}, by $\Gamma$-side-pairing either the fundamental region
or a finite collection of discrete regions congruent to the
fundamental region. Obviously \textit{M} is with
$(\mathbb{U}^3,\Gamma)$-structure.

We shall continue using Hamilton's quaternion $\mathcal{H}$, and the
notation for points \textit{p} in $\mathbb{U}^3$ will be the same as
that in $\mathbb{R}^3$, only with $r>0$.

Furthermore, since we restrict attention to the upper half-space,
the automorphism group becomes $\textrm{PSL}(2,\mathbb{C})$ (linear
fractional transformation with complex coefficients). If \textit{T}
is a map of the form (\ref{mobius transformation}), where
$a,b,c,d\in \mathbb{C}$, we have

\begin{eqnarray}
T(p) &=& T(z+rj) \nonumber \\
     &=&
     \frac{(az+b)(\bar{c}\bar{z}+\bar{d})+a\bar{c}r^2}{\|cp+d\|^2}
     \nonumber \\
     & & +j\frac{r}{\|cp+d\|^2}~.
\end{eqnarray}

For an element $g\in \textrm{PSL}(2,\mathbb{C})$, $g\ne \pm I$ is
classified as follows:

\begin{flushleft}
\begin{enumerate}
  \item[i)] if $|tr(g)|=2$ ~\&~ $tr(g) \in \mathbb{R}$, \textit{T} is
\emph{parabolic};
  \item[ii)] if $|tr(g)|>2$ ~\&~ $tr(g) \in \mathbb{R}$, \textit{T} is
\emph{hyperbolic};
  \item[iii)] if $0\le |tr(g)|<2$ ~\&~ $tr(g) \in \mathbb{R}$, \textit{T} is
\emph{elliptic};
\end{enumerate}
\end{flushleft}

To define explicit expressions for dynamical systems on \textit{M},
we first need to find the so-called automorphic functions that are
invariant under the elements of $\Gamma$.

By definition, an automorphic function \textit{A} for the
\textit{Fuchsian} group $\Gamma$ is a meromorphic function generated
on $\mathbb{U}^3$ such that

\[
A\big(T_i(p) \big)=A(p)
\]

\noindent for all $T_i \in \Gamma$ and $p \in \mathbb{U}^3$
$(p=z+rj)$.

It would be nice if the dynamics on the 3-manifold \textit{M} can be
defined as

\begin{equation}\label{automorphic function}
\dot{p}=A(p),
\end{equation}

\noindent where \textit{A} is an automorphic function. However,
since we are dealing with vector fields, the solutions generated by
(\ref{automorphic function}) in $\mathbb{U}^3$ are not
$\Gamma$-invariant in the sense that dynamics at the boundary of the
fundamental region won't match up when applying the
$\Gamma$-side-pairing. In order to obtain systems $\dot{p}=f(p)$
with $\Gamma$-invariant trajectories, we require the following
invariance of the vector field \textit{f}:

\begin{lemma}
The system
\[
\dot{p}=f(p)
\]
\noindent will have $\Gamma$-invariant trajectories for any given
discrete group ~$\Gamma$ of isometries of hyperbolic $3$-space
\textit{X}, if
\begin{equation}\label{dynamical system}
f(p)=\frac{d\Big(T^{-1}\big(T(p) \big) \Big)}{dp}\cdot f\big(T(p)
\big), \qquad \forall ~T\in \Gamma.
\end{equation}
\end{lemma}

\begin{figure}[!hbp]
\centering
\includegraphics[width=2in]{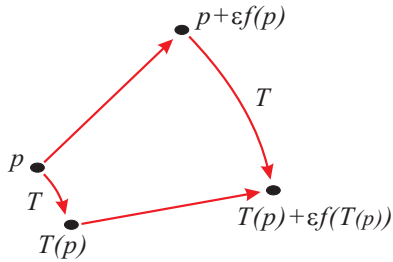}
\caption{Mapping the vector field under any $\textit{T}\in
\Gamma$}\label{matching vector fields}
\end{figure}

\noindent \textbf{Proof.} To make the dynamics match up after the
side-pairing, we require the ``ends'' of infinitesimal vectors in
the direction of $f(p)$ to map appropriately under $\Gamma$ (see
fig.~\ref{matching vector fields} for illustration).

Hence we require

\[
T\big(p+\varepsilon f(p) \big) = T(p) + \varepsilon f\big(T(p) \big)
\]

\noindent for sufficiently small $\varepsilon$. Thus

{\setlength\arraycolsep{1.5pt}
\begin{eqnarray*}
f(p) &=& \frac{T^{-1}\Big(T(p)+\varepsilon f\big(T(p) \big)
\Big)-p}{\varepsilon} \\
f(p) &=& \frac{T^{-1}\Big(T(p)+\varepsilon f\big(T(p) \big) \Big) -
T^{-1}\big(T(p) \big)}{\varepsilon f\big(T(p) \big)} \\
& & \cdot f\big(T(p)\big)
\end{eqnarray*}}

\[
\therefore \qquad f(p)=\frac{dT^{-1}}{dp}\big(T(p) \big) \cdot
f\big(T(p) \big),
\]

\noindent so the lemma is proved. \qquad $\square$

To work out the relation between $f(p)$ and $f\big(T(p) \big)$
explicitly, we have, from (\ref{mobius transformation}),

\[
T^{-1}(p)=\frac{dp-b}{-cp+a}
\]

\[
\therefore \qquad
\frac{d}{dp}\big(T^{-1}(p)\big)=\frac{ad-bc}{(a-cp)^2}
\]

\[
\Longrightarrow \qquad \frac{dT^{-1}}{dp}\big(T(p) \big) =
\frac{(cp+d)^2}{ad-bc}
\]

Therefore, for such a map $\textit{T}\in \Gamma$, the invariance of
the dynamical system \textit{f} given by (\ref{dynamical system})
can be written in the form

\begin{equation}\label{modified dynamical system}
F\big(T(p) \big)=\frac{ad-bc}{(cp+d)^2}\cdot F(p)
\end{equation}

\noindent Note that (\ref{modified dynamical system}) differs from
the scalar invariance

\[
A\big(T(p) \big)=A(p), \qquad \textit{T}\in \Gamma,
\]

\noindent which is given by any automorphic function. So we shall
obtain vector fields \textit{F} that satisfies (\ref{modified
dynamical system}) by modifying the \textit{Poincar\'e} theta series
(see [Ford, 1929]) which can be used to generate automorphic
functions for those \textit{Fuchsian} groups with infinite elements.

\begin{definition}
Let \textit{H} be a rational function, which has no poles at the
limit points of the isometry group $\Gamma$, the theta series is
given by

\begin{displaymath}
\theta(p) = \sum_{i=0}^{\infty} (c_ip+d_i)^2 \cdot H(p_i),
\end{displaymath}

\noindent where $p\in \mathbb{U}^3$, $I,T_1,T_2,T_3,\cdots$ are the
elements of ~$\Gamma$, and

\[
p_i = T_i(p) = \frac{a_ip+b_i}{c_ip+d_i}.
\]

\end{definition}

It is easy to verify that

\[
\theta(p_i) = (c_ip+d_i)^{2m} \cdot \theta(p)
\]

\noindent for each \textit{i}, and by definition, two distinct theta
series $\theta_1$ and $\theta_2$ with the same choice on \textit{m},
we can have

\[
F(p) = \frac{\theta_1(p)}{\theta_2(p)}.
\]

Moreover,

\[
F(p_i)=F(p)
\]

\noindent for each \textit{i}, i.e., \textit{F} is an automorphic
function.

From (\ref{modified dynamical system}), we know that in the case of
dynamical systems, some modification must be made to the theta
series so that they can provide the invariance of the vector fields.
Therefore instead of ~$\theta_1$, we define

\[
\widetilde{\theta}_1(p) =
\sum_{i=0}^{\infty}\frac{(c_ip+d_i)^{2-2m}}{(a_id_i-b_ic_i)}\cdot
H_1\big(T_i(p) \big)
\]

\noindent while keep ~$\theta_2(p)$ as usual.

\begin{lemma}
The function

\[
F(p) = \frac{{\widetilde{\theta}}_1(p)}{\theta_2(p)}
\]

\noindent satisfies

\[
F\big(T_i(p) \big) = \frac{a_id_i-b_ic_i}{(c_ip+d_i)^2} \cdot F(p)
\]

\noindent for each \textit{i} and so defines a $\Gamma$-invariant
dynamical system if $m \ge 2$.
\end{lemma}

\noindent \textbf{Proof.} Since $\theta_2$ is the normal theta
series, we have

\[
\theta_2\big(T_j(p) \big) = (c_jp+d_j)^{2m} \cdot \theta_2(p)
\]

\noindent for each \textit{j}, while for ~$\widetilde{\theta}_1$, we
have

{\setlength\arraycolsep{1.5pt}
\begin{eqnarray*}
\widetilde{\theta}_1\big(T_j(p) \big) &=& \sum_{i=0}^{\infty}\frac{H_1\big(T_iT_j(p) \big)}{\Big(c_i (\frac{a_jp+b_j}{c_jp+d_j})+d_i \Big)^{2m-2}(a_id_i-b_ic_i)} \\
&=& \sum_{i=0}^{\infty}
\frac{(c_jp+d_j)^{2m-2}}{\big((c_ia_j+d_ic_j)p+c_ib_j+d_id_j
\big)^{2m-2}} \\
& & \cdot \frac{H_1\big(T_iT_j(p) \big)}{(a_id_i-b_ic_i)} \\
&=& (c_jp+d_j)^{2m-2}(a_jd_j-b_jc_j)\cdot \\
& & \sum_{i=0}^{\infty} \frac{1}{\big((c_ia_j+d_ic_j)p+c_ib_j+d_id_j
\big)^{2m-2}} \\
& & \cdot \frac{H_1\big(T_iT_j(p)
\big)}{(a_id_i-b_ic_i)(a_jd_j-b_jc_j)}
\\
&=& (c_jp+d_j)^{2m-2}(a_jd_j-b_jc_j)\cdot \widetilde{\theta}_1(p),
\end{eqnarray*}}

\noindent since

{\setlength\arraycolsep{1.5pt}
\begin{eqnarray*}
T_iT_j(p) &=& \frac{a_i \frac{a_jp+b_j}{c_jp+d_j}+b_i}{c_i
\frac{a_jp+b_j}{c_jp+d_j}+d_i} \\
&=&
\frac{(a_ia_j+b_ic_j)p+(a_ib_j+b_id_j)}{(c_ia_j+d_ic_j)p+(c_ib_j+d_id_j)},
\end{eqnarray*}}

\noindent and

{\setlength\arraycolsep{1.5pt}
\begin{eqnarray*}
\textrm{det}\big(T_iT_j(p) \big) &=& (a_ia_j+b_ic_j)(c_ib_j+d_id_j)
\\ & & -(a_ib_j+b_id_j)(c_ia_j+d_ic_j) \\
&=& (a_jd_j-b_jc_j) \cdot (a_id_i-b_ic_i).
\end{eqnarray*}}

\noindent Hence

{\setlength\arraycolsep{1.5pt}
\begin{eqnarray*}
F(p) &=& \frac{\widetilde{\theta}_1(p)}{\theta_2(p)} \\
&=& \frac{\frac{1}{(c_jp+d_j)^{2m-2}(a_jd_j-b_jc_j)} \cdot
\widetilde{\theta}_1\big(T_j(p) \big)}{\frac{1}{(c_jp+d_j)^{2m}}
\cdot \theta_2\big(T_j(p) \big)} \\
&=& \frac{(c_jp+d_j)^2}{a_jd_j-b_jc_j} \cdot F\big(T_j(p) \big) \\
&=& \big((T_j)^{-1} \big)' \big(T_j(p) \big) \cdot F\big(T_j(p)
\big),
\end{eqnarray*}}

\noindent therefore the result follows. \qquad $\Box$

\begin{definition}
An \textbf{automorphic vector field} on $\mathbb{U}^3$ is a
meromorphic, hypercomplex valued function \textit{F}, such that it
satisfies (\ref{modified dynamical system}) for each isometry
\textit{T} in the \textit{Fuchsian} group $\Gamma$.
\end{definition}

From the discussion above, we know that such functions \textit{F},
generate dynamics situated on hyperbolic 3-manifolds, which is
written in the form

\[
\dot{p}=F(p).
\]

\noindent The trajectories are $\Gamma$-invariant on any fundamental
region, we can then either ``wrap up'' one of them or choose a
finite number and apply the $\Gamma$-side-pairing, both of which
will give rise to systems sit on the resulting hyperbolic 3-manifold
explicitly.\\

\begin{figure}[!h]
  \centering
  \includegraphics[width=2in]{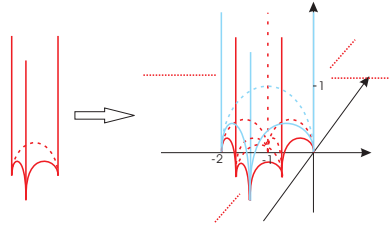}\\
  \caption{Tessellation of $\mathbb{U}^3$ by hyperbolic tetrahedra}\label{tessellation}
\end{figure}

\noindent \textbf{Example.} It is known that the upper half-space
$\mathbb{U}^3$ can be tessellated by hyperbolic ideal tetrahedron.
Fig.~\ref{tessellation} shows one particular representation.

Let the $\Gamma$-side-pairing be either translations or simple
expansions and contractions. According to fig.\ref{tessellation} we
then have the \textit{Fuchsian} group generated by the
transformations

\begin{displaymath}
\begin{array}{ll}
T_1(p)=\frac{p-2}{2}; & \qquad T_2=\frac{p}{2};\\
T_3(p)=\frac{p-1-\sqrt{3}i}{2}; & \qquad
T_4(p)=\frac{p+3+{\sqrt{3}}i}{2}; \\
T_5(p)=p+2.
\end{array}
\end{displaymath}

Choosing

\begin{displaymath}
H_1(p)=p+\frac{1}{2}+\frac{\sqrt{3}}{2}i+5j, \qquad H_2(p)=1.
\end{displaymath}

\begin{figure}[!h]
  \centering
  \includegraphics[width=1.2in, height=1.2in]{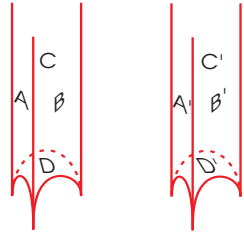}\\
  \caption{Side-pairing two tetrahedra by translation}\label{identification}
\end{figure}

\begin{figure}[!h]
\centering
\includegraphics[width=1.6in,height=1.6in]{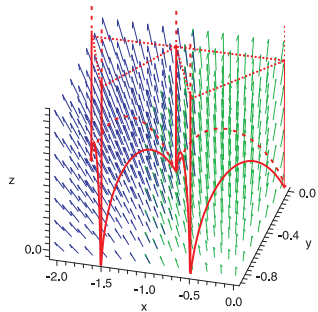}
\caption{The solution trajectories for the system $\dot{p}=F(p)$,
where \textit{F} is generated by $H_1$ and $H_2$.} \label{dynamics}
\end{figure}

We can obtain a dynamical system by using the modified automorphic
functions. Note that in this example, $H_1$ and $H_2$ don't define
poles within the phase-space, however, the system will have poles
introduced by the modified theta series. In fact, the whole
\textit{z}--plane will be covered with equilibria due to the fact
that it contains only cusp points. Fig.~\ref{identification} shows
one possible construction of a hyperbolic 3-manifold by translation.
Moreover, fig.~\ref{dynamics}  illustrates the solution trajectories
of the system (computed in MAPLE), and the vector fields match up
perfectly at the boundaries.

\section{Gluing 3-Manifolds Using the Conformal Ball Model}

We now propose another way of generating dynamical systems on
3-manifolds. Instead of using the upper half-space model, we shall
now investigate hyperbolic 3-manifolds under the conformal ball
model. The same argument applies here, i.e., given a group $\Gamma$
of isometries of $X$ and a proper $\Gamma$-side-pairing, we can form
a 3-manifold $M$ with an $(X,\Gamma)$-structure by gluing a finite
number of disjoint convex polyhedra. Moreover, if we take into
consideration of the dynamical systems naturally situated on those
solid fundamental polyhedra, the $\Gamma$-side-pairing will then
yield a new system defined on the resulting manifold $M$ if and only
if the trajectories match up according to the gluing pattern.

\begin{figure}[!h]
  \centering
  \includegraphics{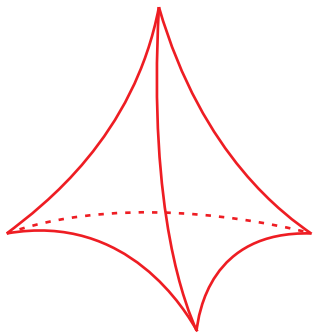}\\
  \caption{An ideal tetrahedron}\label{ideal tetrahedron}
\end{figure}

Again, as an example, we consider a regular ideal tetrahedron in
$B^3$, which has the shape in fig.~\ref{ideal tetrahedron}.

\begin{figure}[!h]
  \centering
  \includegraphics{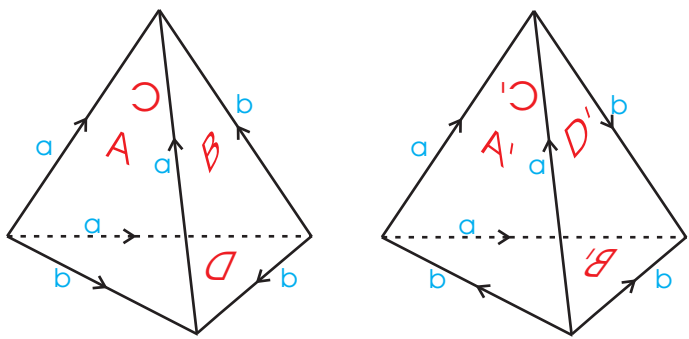}\\
  \caption{The gluing pattern of two regular ideal tetrahedrons}\label{two tetrahedrons}
\end{figure}

Let $T_1$ and $T_2$ be two disjoint regular ideal tetrahedrons in
$B^3$, illustrated in fig.~\ref{two tetrahedrons}. For
simplification, we regard them as regular tetrahedrons in the
Euclidean space.

Because a \textit{M\"obius} transformation of the unit ball $B^3$
leaves it invariant, the permutation of the four vertices will
determine the gluing pattern accordingly. If we label the sides and
edges of $T_1$ and $T_2$ as in fig.~\ref{two tetrahedrons}, there
must exist an isometry of $B^3$ such that the sides of $T_2$,
namely, $A'$, $B'$, $C'$, $D'$, are mapped onto those of $T_1$,
i.e., $A$, $B$, $C$, $D$, and exactly in this order. It can be
proved that this side-pairing is proper, hence implies that the
resulting space will be a hyperbolic 3-manifold, say $M$, which is
known as the figure-eight knot complement.

\begin{figure}[!h]
  \centering
  \includegraphics{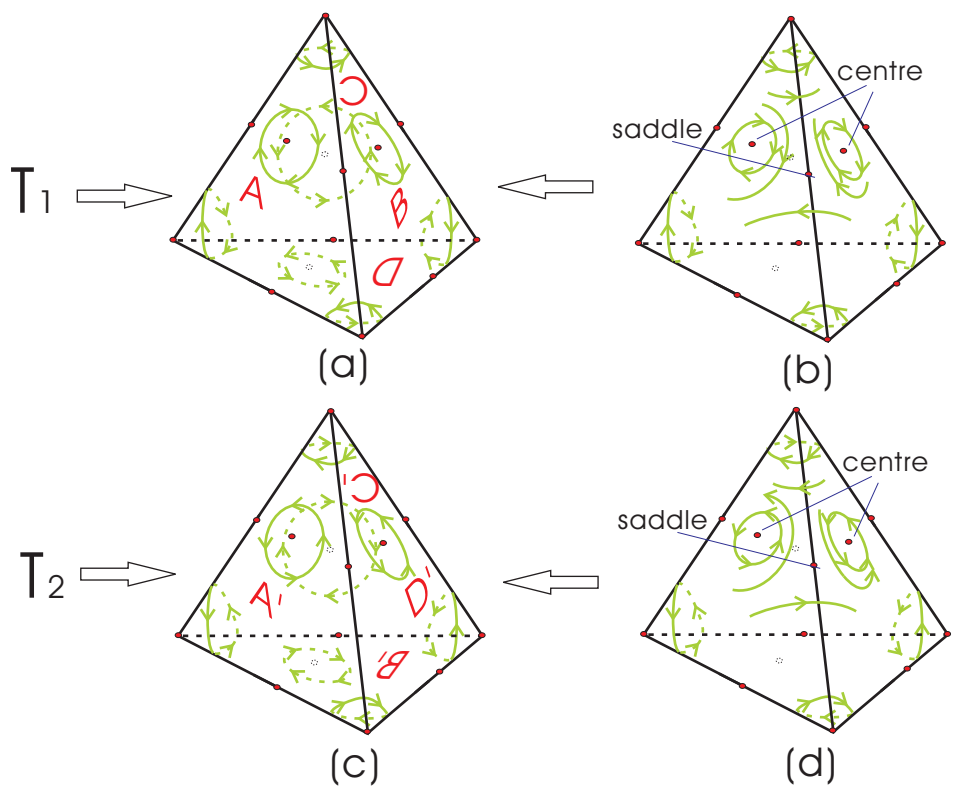}\\
  \caption{Dynamical systems on T1 and T2}\label{two systems}
\end{figure}

Now by assuming the existence of systems on these solid regular
tetrahedrons, a new dynamical system can then be constructed on the
resulting manifold via the side-pairing if and only if the
trajectories match up on the corresponding boundaries of the
polyhedra components. As an example, fig.~\ref{two systems}
illustrates this matching up by applying the side-pairing that we
mentioned above. Note that the explicit dynamics in (a) and (c) are
obtained by repeating (b) and (d) on all sides and edges of $T_1$
and $T_2$, respectively.

\section{Modified Reeb Foliations and Systems on 3-Manifolds}

The classical Reeb foliation of the sphere and the torus are
well-known (see [Moerdijk \& Mrcun, Candel \& Conlon, 2000]). These
are obtained first from a Heegaard splitting of the sphere

\[
S^3\cong X \cup_{\partial X} X
\]

\noindent where $X$ is a solid torus and each copy of $X$ carries
the foliation shown below in fig.~\ref{reeb foliation}.

\begin{figure}[!h]
\centering
\includegraphics{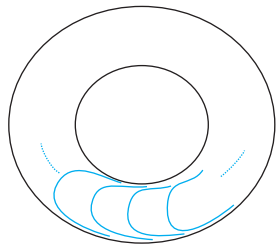}
\caption{The Reeb Foliation.} \label{reeb foliation}
\end{figure}

Each leaf apart from the bounding torus is a plane immersed into the
solid torus. In this paper we shall show that an infinite set of
dynamical systems exists on the 3-sphere which are formed by taking
a genus $p$ (for any $p\ge 1$) Heegaard Splitting of $S^3$ and
finding a generalized Reeb foliation on the solid genus $p$ bounded
3-manifolds. Each leaf (apart from the bounding genus $p$ surfaces
and a singular leaf) will be an unbounded surface of infinite genus.
Of course, it is well-known that every compact three-manifold has a
(nonsingular) foliation (see [Candel \& Conlon, 2000]), essentially
proved by Dehn surgery on embedded tori, each of which carries a
Reeb component. However, this is an existence result and it is
difficult to use to define explicit dynamical systems on
three-manifolds.

\begin{figure}[!h]
\centering
\includegraphics[height=1.2in]{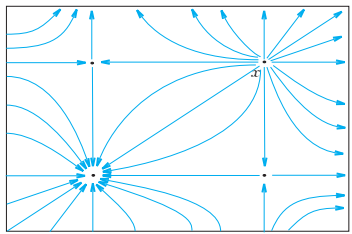}
\caption{A Simple system on the torus} \label{torus system}
\end{figure}

We begin by describing a simple system on the torus which can be
mapped onto each leaf of the Reeb foliation to give a system on
$\mathbb{R}^2$ with an infinite number of equilibria. The basic
system on the torus will consist of a source, a sink and two saddles
as shown in fig.~\ref{torus system}.

\begin{figure}[!h]
\centering
\includegraphics{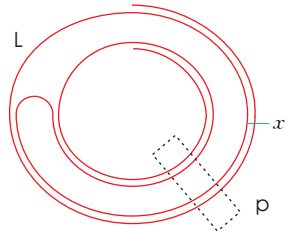}
\caption{A single leaf $L$} \label{a single leaf}
\end{figure}

\noindent (Note that the converse of the \textit{Poincar$\acute{e}$}
index theorem is not true, so it is not possible to have just a
source and a saddle on the torus, although their total index would
be $0$.) Consider a single noncompact leaf in the Reeb foliation
consisting of a `rolled up' plane as in fig.~\ref{a single leaf},

\begin{figure}[!h]
\centering
\includegraphics[width=1.75in,height=1.75in]{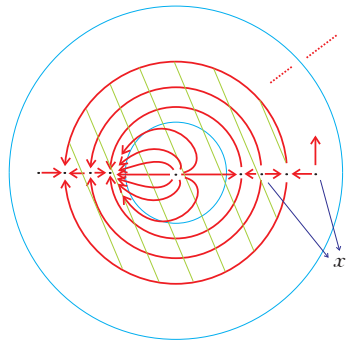}
\caption{Resulting dynamics on the cylinder}\label{source}
\end{figure}

\noindent The plane $P$ cuts the leaf $L$ into an infinite number of
cylinders plus a disk. Mapping the dynamics of fig.~\ref{torus
system} onto each cylinder and adding a source at the origin of the
disk gives the system on the plane shown in fig.~\ref{source}.

We shall organize the dynamics on the leaf so that the sources lie
`below' the point $x$ on the torus when the leaf is folded up.

\begin{figure}[!h]
\centering
\includegraphics{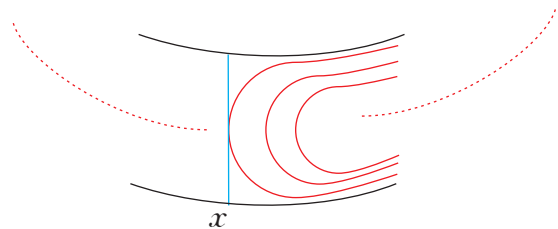}
\caption{Shrinking of the leaf} \label{shaded region}
\end{figure}

Note that the size of the shaded region in fig.~\ref{source} depends
on the leaf and shrinks to zero with origin `below' $x$ as in
fig.~\ref{shaded region}.

We shall now show that there is a (singular) foliation of a
3-manifold of genus $p$ containing a compact leaf consisting of the
bounding genus $p$ surface, an uncountable number of unbounded
leaves of infinite genus and a set of one-dimensional singular
leaves. Consider first the genus $2$ case.

\begin{lemma}
Consider the orientable $3$-manifold with boundary consisting of the
closed surface of genus $2$. There is a singular foliation of this
manifold defined by a dynamical system with a singular
one-dimensional invariant submanifold, an infinite number of
noncompact invariant submanifolds of infinite genus and a single
leaf consisting of the boundary. \label{reeb-like system}
\end{lemma}

\begin{figure}[!h]
\centering
\includegraphics{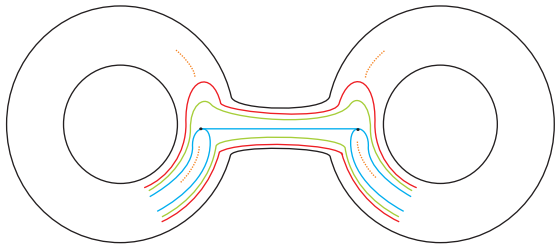}
\caption{Gluing two tori via the leaves} \label{gluing pattern}
\end{figure}

\noindent \textbf{Proof.} We obtain the foliation by modifying the
Reeb foliation and its associated dynamical system introduced above.
Hence consider two systems of the form in fig.~\ref{source}, where
one has the arrows reversed (i.e. we reverse time in the
corresponding dynamical system). We then form the connected sum of
the bounding tori by removing a disk around the source (or sink) at
the point $x$. Then we `plumb' each leaf in a similar way (again
removing the source or sink). This will require one singular line
joining the origins of the leaves which occur just `below' $x$. See
fig.~\ref{gluing pattern} for illustration.

\noindent The leaves clearly have the form stated in the lemma.
\qquad $\Box$
\\

\begin{figure}[!h]
\centering
\includegraphics[width=3in,height=1.5in]{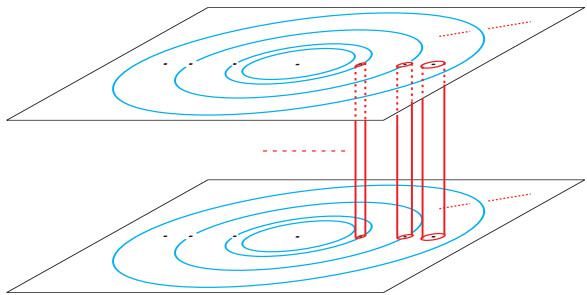}
\caption{A typical leaf} \label{typical leaf}
\end{figure}

\noindent \textit{Remarks.} The nonsingular leaves (apart from the
genus $2$ boundary surface) are embeddings of the surfaces shown in
fig.~\ref{typical leaf}.

Note that we must have at least one singular fibre in order to
introduce such a foliation on a higher genus surface. For we have

\begin{theorem}
Any foliation of codimension $1$ of a compact orientable manifold
$M$ of dimension $3$ with finite fundamental groups and genus $>1$,
which is transversally oriented, must have a singular leaf.
\end{theorem}

\noindent \textbf{Proof.} By Novikov's theorem (see [Moerdijk \&
Mrcun, 2003]), any codimension $1$ transversely orientable foliation
of $M$ has a compact leaf and if $M$ is orientable, this compact
leaf is a torus containing a Reeb component. Thus, if $M$ contains a
compact leaf of genus $>1$, it is not a torus and hence there must
exist a singular leaf. \qquad $\Box$

\noindent \textit{Remarks.} We can find a similar singular foliation
of a genus $2$ 3-manifold by adding a handle between the stable and
unstable points on the torus in fig.~\ref{torus system}. This gives
a typical leaf shown in fig.~\ref{adding handles}, rather than the
one in fig.~\ref{typical leaf}.

\begin{figure}[!h]
\centering
\includegraphics{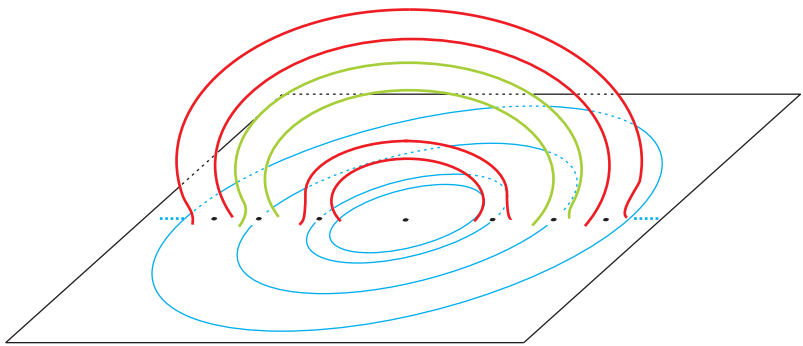}
\caption{A typical leaf obtained by adding handles} \label{adding
handles}
\end{figure}

We now define systems on 3-manifolds by gluing two systems of the
form above situated on solid genus-$p$ surfaces by the use of a
Heegaard diagram. We first recall the general theory of Heegaard
Splittings of 3-manifolds (see, e.g. [Hempel, 1976]). A
\textit{Heegaard Splitting} of genus $p$ $(V_1,V_2)$ of a 3-manifold
$M$ is a pair of solid cubes with $p$ handles $V_1$, $V_2$ such that
$M$ is obtained from $V_1$ and $V_2$ by gluing $\partial V_1$ to
$\partial V_2$. Using a simplicial decomposition of $M$ and a dual
complex, it can be seen that any 3-manifold has a Heegaard
Splitting. Let $\{ D_1,D_2,\cdots, D_n \}$ be pairwise disjoint
properly embedded $2$ cells in $V_2$ which cut $V_2$ into a 3-cell.
Then $\{
\partial D_1,\partial D_2,\cdots, \partial D_n \}$ cut $\partial V_2 = \partial
V_1$ into a 2-sphere with $2n$ holes. We call $(V_1; \partial D_1,
\cdots, \partial D_n)$ a \textit{Heegaard diagram} of $(V_1,V_2)$.
We can get back to $M$ from a Heegaard diagram in the following way:

\begin{enumerate}
\item[(i)] Attach a copy of $B^2\times I$ to $V_1$ ($B^2$ is the 2-ball,
$I=[0,1]$) for each $i=1,\cdots, n$ by identifying $\partial B^2
\times I$ with a neighbourhood of $\partial D_i$ in $\partial V_1$.
The resulting manifold $M_1$ has a 2-sphere boundary.
\item[(ii)] Attach a copy of $B^3$ ($=$3--ball) to $M_1$ via $\partial
B^3$ to $\partial M_1$. This gives $M$.
\end{enumerate}

We can now state

\begin{figure}[!h]
\centering
\includegraphics[height=2in]{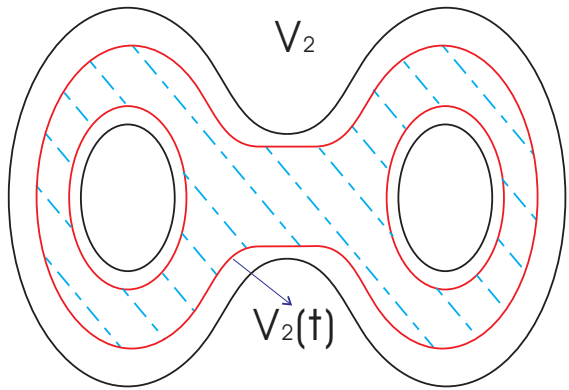}
\caption{A solid genus-$p$ handle-body contained within $V_2$}
\label{solid handle-body}
\end{figure}

\begin{theorem}
For any $3$-manifold $M$, and any $p>0$, there is a Reeb-like
dynamical system on $M$ given by gluing two systems of the form
given in Lemma $\ref{reeb-like system}$.
\end{theorem}
\noindent \textbf{Proof.} Let $(V_1,V_2)$ be a Heegaard Splitting of
$M$ of genus $p$ and let $\phi_1, \phi_2$ be dynamical systems
defined on $V_1,V_2$, respectively, of the form given in Lemma
\ref{reeb-like system}. Let $\psi: \partial V_1 \to
\partial V_2 \simeq \partial V_1$ be the homeomorphism defined in
(i), (ii) above. By using C-homeomorphisms of the type in
[Lickorish, 1962], we can assume that $\psi$ is smooth. Now let
$V_2(t)$ be a solid genus-$p$ handle-body contained within $V_2$ (as
in fig.~\ref{solid handle-body}) so that $V_2(1)=V_2$ and $V_2(0)$
is a solid genus-$p$ handle-body properly contained in $V_2$. We can
extend $\psi$ to a smooth map $\widetilde{\psi}:V_2\to V_2$ by the
homotopy

\begin{equation}
\widetilde{\psi}= \left\{
\begin{array}{ll}
(1-t)I+t\psi & \textrm{on} ~\partial V_2(t) \\
I & \textrm{on}~ V_2(0)
\end{array} \right.
\label{homotopy}
\end{equation}

Let $X_2$ be the vector field corresponding to $\phi_2$ on $\partial
V_2$. Then we `twist' the dynamics on $V_2$ by $\psi$, i.e.,
$(\psi^{-1})_\ast X_2$ and extend this to $V_2$ in an obvious way
using (\ref{homotopy}). Then the dynamics on $\partial V_2$ match
those on $\partial V_1$ according to the Heegaard diagram and the
result is proved. \qquad $\Box$

\section{Conclusions}

In this paper, we have considered a variety of methods for
generating systems on 3-manifolds. We have shown how to construct
dynamical systems explicitly on hyperbolic 3-manifolds. This is
achieved by using a modified theta series to obtain the
`generalized' automorphic functions which `uniformize' the vector
fields on the manifold. Here we concentrated on using the upper
half-space model for the hyperbolic space, while it is also possible
to use the disk model. Also we gave an example of how to generate
such systems. Also we consider constructing dynamical systems with
the help of Reeb foliation. This is achieved by defining a system on
each leave and then using the connected sum method to link them
together.

In the next paper we shall consider the possible existence of
knotted chaotic systems when applying the side-pairing to obtain the
3-manifolds.

\end{document}